\documentclass[10pt,oneside,a4paper,leqno,final]{article}

\usepackage{amsmath}
\usepackage{amsthm}
\usepackage{amssymb}
\usepackage{mathrsfs}
\usepackage{subfigure}
\usepackage{graphicx}

\usepackage{multirow,rotating}

\input xy
\xyoption{all}
\usepackage{geometry}

\newtheorem{thm}{\sc Theorem}[section]      
\newtheorem{cor}[thm]{\sc Corollary}        
\newtheorem{lem}[thm]{\sc Lemma}            
\newtheorem{prop}[thm]{\sc  Proposition}     
\newtheorem{defn}[thm]{\sc Definition}      
\newtheorem{rem}[thm]{\sc Remark}       

\newcommand{\R}{\mathbb R}

\usepackage{bm}
\newcommand{\simbolovettore}[1]{{\boldsymbol{#1}}}
\newcommand{\va}{\simbolovettore{a}}

\newcommand{\vb}{\simbolovettore{b}}
\newcommand{\vp}{\simbolovettore{p}}
\newcommand{\vP}{\simbolovettore{P}}
\newcommand{\vq}{\simbolovettore{q}}

\newcommand{\vx}{\simbolovettore{x}}
\newcommand{\vy}{\simbolovettore{y}}
\newcommand{\vz}{\simbolovettore{z}}
\newcommand{\vg}{\simbolovettore{g}}
\newcommand{\vG}{\simbolovettore{G}}
\newcommand{\vh}{\simbolovettore{h}}
\newcommand{\vi}{\simbolovettore{i}}

\newcommand{\vs}{\simbolovettore{s}}
\newcommand{\vS}{\simbolovettore{S}}

\newcommand{\vj}{\simbolovettore{j}}
\newcommand{\vk}{\simbolovettore{k}}
\newcommand{\vl}{\simbolovettore{l}}
\newcommand{\vf}{\simbolovettore{f}}
\newcommand{\vc}{\simbolovettore{c}}
\newcommand{\vZ}{\simbolovettore{Z}}

\newcommand{\zero}{\boldsymbol{0}}
\newcommand{\Ddt}{\mathrm{\dfrac{ D}{dt}}}
\newcommand{\Ddtp}{\mathrm{\tfrac{ D}{dt}}}

\newcommand{\finedim}{\hfill $\Box$}

\begin{document}
\pagenumbering{arabic}

\title{
  Global dynamics under a weak potential on a sphere \thanks{Work
    partially supported by``Progetto 5 per mille per la ricerca''
    (Bando 2011). ``Collisioni fra vortici puntiformi e fra filamenti
    di vorticita': singolarita', trasporto e caos.''}}

\author{Roberto Castelli,  Francesco Paparella and Alessandro Portaluri}
\author{Roberto Castelli \thanks{BCAM - Basque Center for Applied
    Mathematics, Bizkaia Technology Park, 48160 Derio, Bizkaia, Spain
    ({\tt rcastelli@bcamath.org}).}  \and \and Francesco Paparella
  \thanks {Dipartimento di Matematica ``Ennio De Giorgi'', Ex-collegio
    Fiorini, University of Salento, 73100 Lecce, Italy ({\tt
      francesco.paparella@unisalento.it}).  } \and Alessandro
  Portaluri\thanks{Dipartimento di Matematica ``Ennio De Giorgi'',
    Ex-collegio Fiorini, University of Salento, 73100 Lecce, Italy
    ({\tt alessandro.portaluri@unisalento.it}).}  }

\maketitle

\begin{abstract}
  We give a detailed analytical description of the global
  dynamics of a point mass moving on a sphere under the action of a
  logarithmic potential. After performing a McGehee-type blow-up in
  order to cope with the singularity of the potential, we investigate
  the rest-points of the flow, the invariant (stable and unstable)
  manifolds and we give a complete dynamical description of the
  motion.

\end{abstract}
\noindent {\em MSC Subject Class\/}: Primary 70F10; Secondary 37C80.
\vspace{0.5truecm}


\section*{Introduction} \label{sec:intro}

Topologically, two dimensional Riemann surfaces with constant
(Gaussian) curvature $K$ are classified into three categories:
Euclidean spheres, $\mathbb S^2$ ($K>0$); Euclidean planes, $E^2$
($K=0$); and hyperbolic planes $H^2$ ($K<0$).  Among them, $\mathbb
S^2$ and $E^2$ are more familiar and come out very often in
practice. For example, the mechanics of thin fluid layers on $\mathbb
S^2$ provides a global model of a planetary atmosphere, and on $E^2$
its local approximation.

In this paper we analyze the motion of a point particle moving on a
sphere under the action of a logarithmic potential. Two are the main
reasons for the choice of this particular potential.  First, it arises
in different physical scenarios: such as in models of astrodynamics,
\cite{tremaine}, \cite{stoica}; in the dynamics of a charged particle
in a cylindrically symmetric electric field \cite{hooverman} and in
the mathematical theory of vortex filaments of an ideal fluid
\cite{Newton}, \cite{ponce}.  The second reason relies on the fact
that the logarithmic potential $V(x)=-\log(|x|)$ could be considered
as a limit case for $\alpha\rightarrow 0$ of the homogeneous
potentials $V_{\alpha}=|x|^{-\alpha}$ and, while the latter have been
extensively studied by different authors, the former has not been so
deeply investigated.  In particular one could be interested to know if
(and how) some features regarding for instance the regularization of
collisions, the minimality properties of the solutions, the stability
character, may be extended from the homogeneous to the logarithmic
potential case.  Results in this direction have been achieved for
instance in \cite{cas2}, \cite{gencoll}, \cite{tremaine}.

In addition, we consider a sphere, rather than the classical two or
three-dimensional Euclidean space, as the configuration space. Our
goal is to understand which aspects of the dynamics are affected if
the geometry of the underlying space changes, or equally well, what
survives of the planar dynamics if one considers a curved manifold.

From a dynamical point of view the most interesting feature, and the
hardest obstacle for a full understanding of the motion, is played by
the presence of the singularity in the potential function. Indeed, as
it often happens in celestial mechanics, the singularities are the
source of a complicated dynamics and sometimes they are even
responsible for a sort of chaotic motion. From the mathematical point
of view the singularities represent a severe technical hurdle to
overcome and different techniques have been proposed to regularize the
vector field, , mainly for the homogeneous potential case
\cite{levicivita}, \cite{McGehee74}, \cite{easton}, \cite{gronchi},
\cite{cas2} and \cite{cate}.

This paper is inspired by the recent a work \cite{stoica} that studies
the planar motion of a point mass subject to a logarithmic potential
in an astrodynamic context. To overcome the singularity of the vector
field we adapt to our problem the celebrated {\em McGehee
  transformation\/}, a regularizing change of variables currently
popular in the field of Celestial Mechanics and first introduced in
1974 by McGehee\cite{McGehee74} to solve the collisions in the
collinear three-body problem.

The McGehee transformations consist of a polar type change of
coordinates in the configuration space, together with a suitable
rescaling of the momentum. In this way the total collision is blown-up
into an invariant manifold called {\em total collision manifold\/} over
which the flow extends smoothly. Furthermore, each hypersurface of
constant energy has this manifold as a boundary.  By rescaling time in
a suitable way, it is possible to study qualitative properties of
the solutions close to total collision, obtaining a precise
characterization of the singular solutions.

The McGehee transformation are usually applied to the case of
homogeneous potentials but, as shown in \cite{stoica} and as it will be
manifest throughout this paper, with slight modifications they give
interesting results even in the presence of a logarithmic potential.
In fact, although the lack of homogeneity of the logarithmic
nonlinearity breaks down some nice and useful properties of the
transformation, it is still possible to regularize the vector field
and therefore it is still possible to carry out a detailed analytical
description of the rest points, of the invariant manifolds, and of the
heteroclinics on the collision manifold.

The paper is organized as follows: first we introduce some basic
notions about the Hamiltonian formulation of the co-geodesics flow on
a general Riemanian manifold, then in Section \ref{sec:stereo} we
restrict to the case of the sphere and we formulate the equivalent
co-geodesics flow through the stereographic projection. In Section
\ref{sec:description} we introduce the singular logarithmic potential 
and we write the equation of motion.  Section \ref{sec:stabunstab}
deals with the in-deep study of the dynamical system: we regularize
the singularity of the potential with the modified McGehee technique,
and we provide an analysis of the flow on the collision and the zero
velocity manifolds.  Section \ref{sec:global} concerns the global
dynamics and we rephrase the results in terms of the original motion
on the sphere with untransformed coordinates.

\tableofcontents

\section{Preliminaries}

Let $M$ be a Riemannian manifold, namely a smooth $n$-dimensional
manifold $M$ endowed with a metric given by a positive definite
(non-degenerate) symmetric two-form $g$. We denote by $D$ the
associated Levi-Civita connection and by $\Ddtp$ the covariant
derivative of a vector field along a smooth curve $\gamma$. Let $I$ be
an interval on the real line and let $V$ be a smooth function defined
on $I\times M$.
\begin{defn}
  A {\em perturbed geodesic\/} abbreviated as {\em p-geodesic\/} is a
  smooth curve $\gamma\colon I \rightarrow M$ which satisfies the
  differential equation
\begin{equation}\label{eq:pgeodesic}
\Ddt \gamma '(t) +\nabla V(t,\gamma(t))=0
\end{equation}
where  $\nabla V$ denotes the
gradient of $V (t,-) $ with respect to the metric $g.$
\end{defn}
\begin{rem}
  From a dynamical viewpoint, the data $(g,V)$ define a mechanical
  system on the manifold $M$, with kinetic energy $\frac{1}{2}g(v ,v)$
  and time dependent potential energy $V.$ Solutions of the
  differential equation \eqref{eq:pgeodesic} correspond to
  trajectories of particles moving on the Riemannian manifold in the
  presence of the potential $V$. If the potential vanishes we get
  trajectories of free particles and hence geodesics on $M$.  This
  motivates the suggestive name, ``perturbed geodesics" in the case
  $\nabla V\neq 0$. Moreover, if the potential $V$ is time
  independent, modulo reparametrization, perturbed geodesics become
  geodesics of the Jacobi metric associated to $(g,V)$: indeed the
  total energy $$ e = \frac{1}{2}g (\gamma (t)) ( \gamma '(t) , \gamma
  '(t) ) + V (\gamma(t) )$$
  is constant along the any trajectory
  $\gamma$ thus, whenever $V$ is bounded from above, the solutions of
  \eqref{eq:pgeodesic} with energy $e> \sup_{m \in M} V(m)$ are
  nothing but reparametrized geodesics for metric $ [ e-V]g$ on $M$
  with total energy one \cite{AM}.
\end{rem}
Denoting by $(q^1, \dots, q^n)$ a local system of coordinates on $M$,
equation \eqref{eq:pgeodesic} reduces to
\[
 \ddot{q}^i + \Gamma_{jk}^i \dot{q}^j \dot{q}^k = - g^{ij} \dfrac{\partial V}{\partial q^j},
\]
where, as usual, $g^{ij}=(g)^{-1}_{ij}$, and $\Gamma_{jk}^i$ are the
Christoffel symbols.

\subsection*{Geodesic flow as Hamiltonian flow}
The geodesic flow turns out to be a Hamiltonian flow of a special
Hamiltonian vector field defined on the cotangent bundle of the
manifold. The Hamiltonian depends on the metric on the manifold and it
is a quadratic form consisting entirely of the kinetic term.  The
geodesic equation corresponds to a second-order nonlinear ordinary
differential system. Therefore by suitably defining the momenta it can
be re-written as first-order Hamiltonian system.

More explicitly, let us consider a local trivialization chart of the
cotangent bundle $T^*M$
\[
 T^*M \Big\vert_U \cong U \times \R^n
\]
where $U$ is an open subset of the manifold $M$, and the tangent space
is of rank $n$. Let us denote by $(q_1, q_2,\dots , q_n, p_1, p_2,\dots,
p_n)$ the local coordinates 
on $T^{*}M$ and introduce the Hamiltonian
\begin{equation}\label{eq:hamiltgeo}
H: T^*M \to\R:  H(\vq,\vp)= \dfrac12 g^{ij}(\vx)p_i \, p_j\ .
\end{equation}
The Hamilton-Jacobi equations of the geodesic equation with respect to
the metric $\vg$ can be written as
\[
 \left\{
\begin{array}{ll}
\dot q^i = \dfrac{\partial H}{\partial p_i}= g^{ij}(x) p_j\\
\dot p_i = -\dfrac{\partial H}{\partial q_i}= -\dfrac12\dfrac{\partial g^{jk}}{\partial q^i}\,p_j\,p_k\ .
\end{array}\right.
\]
The second order geodesic equations are easily obtained by
substitution of one into the other.  The flow determined by these
equations is called the {\em co-geodesic flow\/}, while the flow
induced by the first equation on the tangent bundle is called {\em
  geodesic flow}.  Thus, the geodesic lines are the projections of
integral curves of the geodesic flow onto the manifold $M$.

Being the Hamiltonian $H$ time-independent, it is readily seen that
the Hamiltonian is constant along the geodesics.  Thus, the
co-geodesic flow splits the cotangent bundle into level sets of
constant energy
\[
 M_E=\{(\vq, \vp) \in T^*M: \quad H(\vq,\vp)=E \},
\]
for each energy $E\geq 0$ , so that
\[
 T^*M=\bigcup_{E\geq0} M_E.
\]
Now let $\vg,\vh$ be two Riemannian metrics on $M$ in the same
conformal class; namely there exists a positive and smooth function
$\lambda=\lambda(\vq)$ of the coordinates such that
\[
\vg^{ij}= \lambda \,\vh^{ij}
\]
or equivalently $\vg=\lambda^{-1}\vh$.  From the definition
\eqref{eq:hamiltgeo} it follows that a scaled co-geodesic Hamiltonian
function corresponds to a conformal change of the metric. In fact
if $H_\vg$ and $H_\vh$ denote the Hamiltonian co-geodesic functions
and if $\vg, \vh$ are in the same conformal class then it immediately
follows by \eqref{eq:hamiltgeo} that
\[
 H_\vg(\vq,\vp)= \lambda\, H_\vh(\vq, \vp).
\]
As a consequence, Hamilton's equations with respect to this two
Hamiltonian functions are related as follows
\begin{equation}\label{eq:conformal}
\left\{
\begin{array}{ll}
\dot \vq = \hphantom{-}\partial_\vp H_\vg= \lambda\partial_\vp H_\vh + H_\vh \partial_\vp\lambda= \lambda\partial_\vp H_\vh\\
\\
\dot \vp = -\partial_\vq H_\vg= -\lambda\partial_\vq H_\vh - H_\vh \partial_\vq \lambda
\end{array}\right.
\end{equation}
where the last equality in the first equation comes by the fact the
the function $\lambda$ only depends on $\vq$.

In the following we consider a perturbed-geodesic flow on a
sphere, thus it is worth to write down explicitly the free geodesic
flow when the manifold $M$ is a surface of revolution in $\R^{3}$.

Denote with $(x,y)$ the Cartesian coordinates of $\R^2$ and consider
the function $\varphi:U \subset \R^2 \to \R^3$ given by $
\varphi(x,y)=(f(y)\cos x, f(y) \sin x, g(y))$
\[
 U=\{(x,y)\in \R^2: 0\leq x< 2\pi, y_0 < y < y_1\},
\]
where $f$ and $g$ are differentiable functions, with $f'(y)^2+ g'(y)^2
\neq 0$ and $f(y)\neq 0$.  
Thus $\varphi(x,y)$ is an immersion and the
image $\varphi(U)$ is the surface generated by the rotation of the
curve $(0, f(y), g(y))$ around the $z$ axis.\footnote{%
  Here we are considering the Euclidean space equipped with Cartesian
  coordinates whose axis are labeled as $x,y,z$ according to the
  ordering induced by the canonical orthonormal basis of $\R^3$.} The induced Riemannian metric 
$g=(g_{ij})$ in the $(x,y)$
coordinates is given by
\[
 g_{11}= f^2\qquad g_{12}=0 \qquad g_{22}=(f')^2 + (g')^2.
\]
From \eqref{eq:hamiltgeo}, the Hamiltonian function associated to the
geodesic flow is given by
\[
 H(x,y, p_x, p_y)=\dfrac12\left(\dfrac{1}{f^2}p_x^2+ \dfrac{1}{f'^2+g'^2} p_y^2\right)
\]
and the co-geodesics flow reads as
\begin{equation}\label{eq:cogeodesicrot}
\left\{
\begin{array}{ll}
\dot x= \dfrac{1}{f^2} p_x\\
\dot y= \dfrac{1}{f'^2+g'^2} p_y\\
\dot p_x=0\\
\dot p_y= \left[\dfrac{f f'}{f^4}p_x^2+ \dfrac{f'f''+ g'g''}{(f'^2+g'^2)^2}\right]       
\end{array}\right.
\end{equation}
or, equivalently, as
\[
\left\{
\begin{array}{ll}
\ddot x + \dfrac{2f\,f'}{f^2}\dot x\, \dot y=0\\
\ddot y - \dfrac{f\,f'}{f'^2+ g'^2} (\dot x)^2+ \dfrac{f'f''+ g'g''}{f'^2+g'^2}(\dot y)^2=0.
\end{array}\right.
\]

\section{The stereographic projection of the sphere}\label{sec:stereo}

\begin{figure}[ht]
\centering
{
\includegraphics[width=0.500\textwidth]{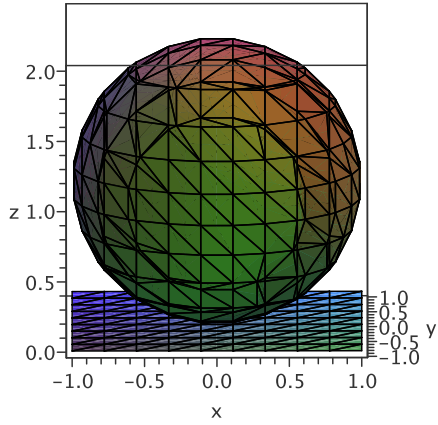}
}
\vspace{-20pt}\caption{Mutual positions of the sphere ($R=1$) and the 0 plane.}
\label{fig:proj}
\end{figure}
It turns out that transformations of the McGehee type may be devised
without too many difficulties for equations which are written in
Cartesian coordinates on a plane. Therefore, rather than attempting to
work directly onto the sphere, we felt it would be more easy (and more
clear) first to project the dynamics on a stereographic plane, and
then to remove the singularities of the resulting equations.

We work on a two-dimensional spherical surface $\vS$ of radius
$R$ and center at the point $C=(0,0, R)$, namely
\[
\vS:=\{(x,y,z) \in \R^3: \ \ x^2+ y^2 + (z-R)^2=R^2\}.
\]
where $(x,y,z)$ are the Cartesian coordinates in $\R^3$ (figure
\ref{fig:proj}).  We shall call {\em north pole\/} and {\em south
  pole\/} the point $N:=(0,\dots, 0,2R) \in \vS$ and its
antipodal $S:=(0, \dots, 0, 0) \in \vS$, respectively. Note that the sphere is
tangent at the origin to the plane $\{\vz=0\}$, that we identify with
$\R^2$.  Next we introduce the stereographic projection
\[\begin{array}{rl}
\pi_\vS: \vS\setminus\{N\} &\longrightarrow \R^2\\
P &\longmapsto \tilde  P, 
\end{array}
\] 
defined by requiring that the three points $N, P, \tilde P$ are
collinear. By a straightforward calculation it follows that the map
$\pi_\vS$ is given explicitly by
\begin{equation}\label{eq:proiezione}
\pi_\vS(x,y,z)=\dfrac{2\,R}{2\,R- z}\,( x, y).
\end{equation}
We use slightly non-standard angular coordinates $\phi$ and
$\theta$ for the spherical surface:
\begin{itemize}
\item $\phi\in[0,2\pi)$ is the usual polar angle of the projection of
  $P$ onto the plane $z=0$;
\item $\theta\in[0,\pi]$ is the angle between the segment
  $\overline{PC}$ and the negative direction of the $z$-axis.
\end{itemize}
A generic point $P=(x,y,z)$ on the sphere in these coordinates has a
local parameterization given by
\[
\left[
\begin{array}{ll}
x\\
y\\
z
\end{array}\right]=R\,  \left[
\begin{array}{ll} 
\sin\theta \cos\phi\\ 
 \sin \theta \sin \phi\\
1-\cos \theta\end{array}\right]=\vP(\phi,\theta).
\]
Of course the map is a diffeomorphism of class $\mathscr C^\infty$. In these 
coordinates the stereographic projection $\pi_\vS$  is defined as :
\[
R  \left[
\begin{array}{ll} 
\sin\theta \cos\phi\\ 
\sin \theta \sin \phi\\
1-\cos \theta\end{array}\right]\longmapsto 
\dfrac{2\, R}{1+\cos\theta}
\left[
\begin{array}{ll} 
\sin\theta \cos\phi\\ 
\sin \theta \sin \phi
\end{array}\right]\ .
\]
We recall that if $M \subset \R^3$ is a portion of a regular surface
represented in Cartesian local coordinates by the vector equation
\[
\vP(u,v):= \zero+ x(u,v)\,\vi+ y(u,v)\,\vj+ z(u,v)\, \vk
\]
then $d\vP= \vP_u\, du + \vP_v\, dv$ for $\vP_u=(x_u, y_u, z_u)$ and
$\vP_v=(x_v, y_v, z_v)$ and hence the metric is given by $ds^2=
\|d\vP\|^2$. With the above 
parametrization of the sphere it follows that
\[
\vP_\theta=R\, (\cos\theta\cos\phi, \cos \theta \sin \phi, \sin\theta), \qquad \vP_\phi=R (-\sin\theta \sin\phi, \sin\theta \cos\phi,0). 
\]
Denoting by $\vg$ and $\vg_\vS$ respectively the Riemannian metric on
the sphere $\vS$ and the metric on the plane induced by the
stereographic projection, we have
\[
\vg:= R^{2}\left[\begin{array}{c c}
 \sin^2\theta\, &0\\
 0  &1
 \end{array}\right], \qquad \vg_\vS:=\dfrac{4}{(1+\cos\theta)^2}\, \vg
\]
As a consequence of the above calculation the following result holds.
\begin{lem}
  The manifolds $(\vS, \vg)$ and $(\overline{\R^2}, \vg_\vS)$ are in
  the same conformal class, where $\overline{\R^2}$ denotes the
  Alexandroff compactification of $\R^2$.
\end{lem}

\subsection*{Co-geodesic flows}

Let $\vG$ and $\vG_\vS$ be the matrices corresponding to the inverse
of $\vg$ and $\vg_\vS$ respectively; thus we have
\[
\vG= R^{-2}\,  
\begin{pmatrix}
\sin^{-2}\theta& 0\\
0 & 1 
\end{pmatrix}, \quad \textrm{and}\quad  \vG_\vS= \dfrac{(1+\cos\theta)^2}{4}\, \vG. 
\]
We denote by $\hat \vS=\vS\setminus\{N,S\}$, the sphere minus the
north and the south pole, and by $T^*\hat\vS$ its cotangent bundle
where it is well-defined the Hamiltonian function
\[\begin{array}{rl}
H_{geod}: T^*\hat \vS &\to \R\\
(\vq; \vp)&\mapsto \dfrac12 \langle \vG \,\vp, \vp  \rangle.
\end{array}
\]
Here $\vq:=(\phi, \theta)$ are the positions and $\vp:=(p_\phi,
p_\theta)$ are the momenta. With this 
choice the Hamiltonian function is given by
\[
 H_{geod}(\phi, \theta;p_\phi, p_\theta)= \dfrac{1}{2\,R^2} \left(\dfrac{1}{\sin^2\theta}\,
p_\phi^2+ p_\theta^2\right)
\]
and, as a particular case of \eqref{eq:cogeodesicrot}, the co-geodesic flow  
on the sphere may be written as 
\begin{equation}\label{eq:HJonsphere}
\left\{
\begin{array}{ll}
\dot \phi = \dfrac{1}{R^2\, \sin^2\theta} \, p_\phi\\
\\
\dot \theta= \dfrac{p_\theta}{R^2}\\
\\
\dot p_\phi=0\\
\\
\dot p_\theta=\dfrac{\cos\theta}{R^2\, \sin^3\theta} \,p_\phi^2\\
\\
\end{array}\right. .
\end{equation}
The co-geodesic flow on the plane $\{\vz=0\}$ is equivalent to the
above one through the stereographic projection.  Since the metrics
$\vg$ and $\vg_{\vS}$ are in the same conformal class, the new system
is easily derived using \eqref{eq:conformal} with
$\lambda=\frac{(1+\cos\theta)^{2}}{4}$.

However, on the plane we prefer to use the Cartesian
coordinates $(x,y)$ rather than the angular coordinates $(\phi,\theta)$.
The latter are related to the former by the transformation
 
\[
\phi=\arctan\left( \dfrac{y}{x}\right), \qquad 
\theta= 2\arctan\left(\dfrac{\sqrt{x^2+y^2}}{2R}\right).
\]

Hence, denoting by $\R^2_s$ the plane endowed with the metric
$\vg_\vS$, the Hamiltonian function for the co-geodesic flow
on $(\R^2, \vg_\vS)$
\[
K_{geod}: T^* \R^2_s\to \R: (\vq; \vp)\mapsto \dfrac12 \langle \vG_\vS \,\vp, \vp  \rangle
\]
is explicitly given by 
\[
K_{geod}(x, y;p_x, p_y)= \,\vl (x,y)\, [\va(x,y)\,
p_x^2+ p_y^2].
\]
where $\vq:=(x, y)$ are the positions and $\vp:=(p_x, p_y)$ are the
momenta. In order to derive this expression we set
\[
\va(x,y):= \left[\dfrac{4R^2+x^2+y^2}{4R\sqrt{x^2+y^2}} \right]^2, \qquad \vl(x,y):=\dfrac{8R^2}{(4R^2+ x^2+ y^2)^2}.
\]
and we use the identities
\[
 \sin\theta= \dfrac{4R\sqrt{x^2+y^2}}{4R^2+x^2+y^2},\qquad
\sin\phi= \dfrac{y}{\sqrt{x^2+y^2}}.
\]
and  
\[
\dfrac{(1+\cos\theta)^2}{8R^2}=  \dfrac{8R^2}{(4R^2+ x^2+ y^2)^2}.
\]
Note that $\va(x,y)$ corresponds to the term $(\sin\theta)^{-2}$, 
while $\vl(x,y)$ is just $\frac{1}{2R^{2}}\lambda$. 
To the Hamiltonian $K_{geod}$ is associated the Hamiltonian flow
\begin{equation}\label{eq:HJonplanenotresc}
\left\{
\begin{array}{ll}
\dot x =2\, (\va\,\vl)(x,y)\, p_x\\
\\
\dot y= 2\vl(x,y)\, p_y\\
\\
\dot p_x=-[\partial_x(\va \, \vl) p_x^2+ (\partial_x \vl)\, p_y^2 ]\\
\\
\dot p_y=- [\partial_y(\va \, \vl) p_x^2+ (\partial_y \vl)\, p_y^2 ]
\end{array}\right. 
\end{equation}
where
\[
 (\va\,\vl)(x,y):= \dfrac{1}{2\,(x^2+y^2)}, \quad \partial_x \vl(x,y)=-\dfrac{32 \,R^2x}{(4R^2+x^2+y^2)^3}, \quad 
 \partial_y \vl(x,y)=-\dfrac{32 \,R^2y}{(4R^2+x^2+y^2)^3}.
\]


\section{Position of the problem}\label{sec:description}

We are now in the position to introduce a conservative force field on
the sphere that perturbs, not necessarily by small amounts, the
geodetic dynamics of a free particle governed by the equations discussed
in the previous section.

We place the singularity of the potential at the point $Q=(0,R,R)$. As
a naming convention, we shall often refer to it as the {\it vortex}
point. On $\vS\setminus\{Q\}$ we define the {\em logarithmic
  potential $U$ \/} as
\[
U(P):= -\dfrac{\Gamma}{4\pi}\log (\|\overline{PQ}\|)
\]
where $\|\overline{PQ}\|$ is the three-dimensional Euclidean distance
between $P$ and $Q$, i.e.  $\|\overline{PQ}\|$ is the length of the
chord between two points on the sphere.  We denote by $\vf$ the force
field generated by the potential $U$, that is $\vf(P)=\nabla U(P)$.
At any point $P\neq Q$ it associates a force pointing towards $Q$ if
$\Gamma>0$ or in the opposite direction otherwise, and proportional to
the inverse of the distance $\|\overline{PQ}\|$. Note that, unlike the
planar case, the force field is not tangent to the manifold: at any
point $P$ one could decompose the force vector into two components,
one directed as the normal to the sphere the other tangent to the
sphere. We take the first as balanced by the smooth constraint given
by requiring that the motion happens on the spherical surface,
therefore only the second contributes to the motion.

The Hamiltonian augmented with the potential function
$$H_{mech}: T^* (\hat \vS\setminus\{Q\}) \to \R$$
in the $(\phi, \theta)$-coordinates is
\[
 H_{mech}(\phi,\theta, p_\phi, p_\theta):=H_{geod}(\phi, \theta, p_\phi, p_\theta) + \dfrac{\Gamma}{8\pi} \log(2 R^2(1-\sin \theta \sin \phi)).
\]
where $2 R^2(1-\sin \theta \sin \phi) = \|\overline{PQ}\|^2$.

On $(\R^{2},\vg_{\vS})$ the distance $\|\overline{PQ}\|^2$ becomes
the function
\[
\qquad \vb(x,y):=\left[\dfrac{2R^2[x^2+(y-2R)^2]}{4R^2+x^2+y^2}\right],
\]
therefore we introduce the Hamiltonian 
$$
K_{mech}: T^*(\R^2_s\setminus\{0,V\})\to \R
$$
\[
K_{mech}(x, y, p_x, p_y):=K_{geod}(x,y, p_x, p_y) + \dfrac{\Gamma}{8\pi} \log\vb(x,y)\ .
\]
We observe that $\log\vb \in \mathscr C^\infty(\R^2\backslash\{(0,2R)\})$: 
indeed the point $V=(0,2R)$
corresponds to the stereographic projection of the vortex $Q\in \vS$,
while the origin is a singularity of the metric.

Hamilton's equations associated to $K_{mech}$ can be written as follows:
\begin{equation}\label{eq:hamiltonequ}
\left\{
\begin{array}{ll}
\dot x =2\, (\va\,\vl)(x,y)\, p_x\\
\\
\dot y= 2\vl(x,y)\, p_y\\
\\
\dot p_x=-[\partial_x(\va \, \vl) p_x^2+ (\partial_x \vl)\, p_y^2 + \Gamma/(8\pi) \vb(x,y)^{-1}\,\partial_x \vb(x,y)]\\
\\
\dot p_y=- [\partial_y(\va \, \vl) p_x^2+ (\partial_y \vl)\, p_y^2 + \Gamma/(8\pi) \vb(x,y)^{-1}\,\partial_y \vb(x,y)]
\end{array}\right. \ .
\end{equation}

These equations govern the motion of a particle constrained on a
spherical surface and subject to a force field generated by a
logarithmic potential, as seen on a stereographic plane.
\begin{rem}
  Before we delve into the analysis of system \eqref{eq:hamiltonequ}
  let us remark that when defining the potential function one could
  consider different notions of the distance between two points on a
  sphere. A reasonable choice could be the geodesic distance, that is
  the length of the shortest arc of a great circle passing through two
  points. More precisely, for any couple $\vx,\vy\in \mathbb S^n(r)
  \subset \R^{n+1}$, the geodesic 
  distance $d_\mathbb S(\vx, \vy)$ is given by
\begin{equation}\label{eq:distanzageodetica}
d_\mathbb S(\vx, \vy):= r\, \arccos \dfrac{\langle \vx, \vy\rangle}{r^2},
\end{equation}
where $<\cdot,\cdot>$ is the Euclidean scalar product in $\R^{n+1}$.
In the following only the chord distance will be considered, but we
guess that the local flow, that is the dynamics close and up to the
singularity, should not be different when the geodesic distance is
taken into account. On the other side, we expect the global flow
to be slightly different. However a complete study of the dynamics
with the geodesic distance, the differences and similarity with the
dynamics on the plane and with the one here described, could be
material for future investigations.
\end{rem}

\section{Energy hypersurfaces, regularization and flow}\label{sec:stabunstab}

We begin the analysis of system \eqref{eq:hamiltonequ}
with the description of the topology of the constant-energy
hypersurfaces associated to $K_{mech}$.  For any $h\in\R$ the
hypersurface of constant energy $h$ is given by
\begin{equation}
\begin{array}{rl}
 \widetilde \Sigma_h&:= \{(x,y, p_x, p_y) \in T^*X: K_{mech}(x,y, p_x, p_y)=h\}\\
 &=\left\{(x,y, p_x, p_y) \in T^*X:   \va(x,y) \, p_x^2+ p_y^2 = \dfrac{1}{\vl(x,y)}
\Big(h - \Gamma/(8\pi) \, \log(\vb(x,y))\Big)\right\}.
\end{array}
\end{equation}
where $X:=\R^2\backslash\{0, V\}$ denotes the configuration space and
$T^*X$ the phase space (the cotangent bundle over $X$).

Since the function $\va(x,y)$ is strictly positive, for any value of
$h$ the motion is allowed only in those regions of the configuration
space where the right hand side of the equation in the definition of
$\widetilde \Sigma_h$ is positive (Fig.\ref{fig:energylevels}).  In the Lemma
\ref{thm:zonedinamiche} below, the analysis is performed for
$\Gamma>0$: changing the sign of $\Gamma$ simply switches the allowed region
with the forbidden region.

\begin{figure}[ht]
\centering
\includegraphics[width=0.5\textwidth]{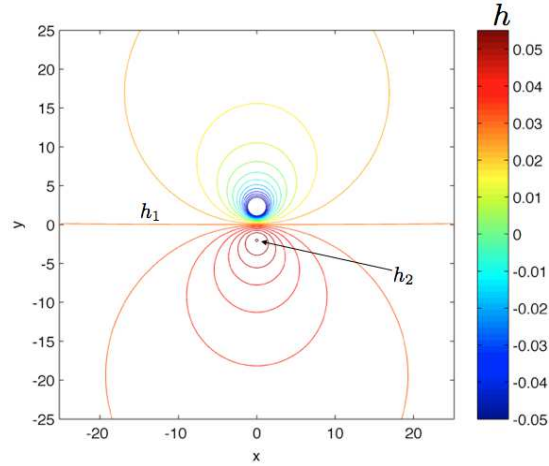}
\caption{Zero-level curves of the function  $\tilde E_{h}(x,y)$ for different values of $h$. ($\Gamma=R=1$)}
\label{fig:energylevels}
\end{figure}

\begin{lem}\label{thm:zonedinamiche}
For any fixed $h$ let $\widetilde E_{h}: \R^{2} \rightarrow \R$ be defined by:
\begin{equation}\label{eq:ehatinalpha}
\widetilde E_{h}(x,y):= \Big(h - \Gamma/(8\pi) \, \log(\vb(x,y))\Big),
\end{equation}
and let $h_{1}=\tfrac{\Gamma}{8\pi}\log(2R^{2})$, $h_{2}=\tfrac{\Gamma}{4\pi}\log(2R)$.
Then 
\begin{enumerate}
\item for every  $h >h_2$ the surface $\widetilde E_{h} (x,y)$ is positive for any  $(x,y)\in\R^{2}$;
\item for any $h \in (h_1, h_2)$ there exists a disk $D^{1}_{h}$ in
  the $\{y<0\}$ half-plane and containing the point $(0,-2R)$ such
  that $\widetilde E_{h} (x,y)$ is positive for each $(x,y) \in
  \R^{2}\backslash D^{1}_h$ and negative otherwise;
\item for any $h < h_1$ there exists a disk $D_h^{2}$ in the $\{y>0\}$
  half-plane, containing the point $(0,2R)$, such that $ \widetilde
  E_{h}(x,y)$ is positive 
  for all $(x,y) \in D_h^{2}$ and negative otherwise.
\end{enumerate}
\end{lem}
\proof The logarithm is a monotone function, therefore the topology of
the level sets of $\widetilde E_{h}(x,y)$ only depends on the level
sets of $\vb(x,y)$.  For any $\delta$ the $\delta$-level set of
$\vb(x,y)$ is given by the point $(x,y)$ lying on the circle
$C_{\delta}$ with equation
$$
x^{2}+y^{2}-\frac{8R^{3}}{2R^{2}-\delta}y+4R^{2}=0\ .
$$
The center of $C_{\delta}$ is placed in the point
$O_{\delta}=(0,\tfrac{4R^{3}}{2R^{2}-\delta})$ and the radius is
$r_{\delta}=\tfrac{2R}{|2R^{2}-\delta|}\sqrt{4\delta
  R^{2}-\delta^{2}}$. It follows that the function $\vb(x,y)$ only
admits values $\delta$ in the range $\delta\in[0,4R^{2}]$. For
$\delta=4R^{2}$ the $\delta$-level set restricts to the point $\tilde
P=(0,-2R)$ while, 
as $\delta$ decrease towards $2R^{2}$, the level set
consists of a circle completely contained in the $\{y<0\}$
half-plane. Moreover it can be checked that
$|O_{\delta}-(0,-2R)|<r_{\delta}$, meaning that the point $\tilde P$
is always surrounded by these circles. The value $\delta=2R^{2}$ is a
singularity for the topology of the level sets: indeed the center of
$C_{\delta}$ as well as the radius $r_{\delta}$ diverge. Note that
$\vb(x,0)=2R^{2}$ for any $x$ and
$\lim_{x^{2}+y^{2}\to\infty}{\vb(x,y)}=2R^{2}$.  Then as $\delta$
decreases below $2R^{2}$ towards zero, the circles $C_{\delta}$ live
in the positive-$y$ halfplane and shrinks around the point $\tilde
Q=(0,2R)$.  \finedim

In order to develop a McGehee type transformation, we define the functions
\begin{equation}\label{eq:lefi}
\left\{
\begin{array}{ll}
\varphi_1(r):= r\, e^{-1/r^2}\\
\varphi_2(r):= 1/r  
\end{array}\right.
\end{equation}
and, following the notation of \cite{stoica}, we introduce the change
of variables
\begin{equation}
\left\{\begin{array}{ll}
x=\varphi_1(r)\, s_1\\
y=  \varphi_1(r)\, s_2 + 2R
\end{array}\right.,\qquad  \left\{\begin{array}{ll}
p_x=\varphi_2(r)\, z_x\\
p_y=  \varphi_2(r)\, z_y 
\end{array}\right.
\end{equation}
where $\vs=(s_1, s_2)=(\cos\alpha,\sin\alpha) \in \mathbb S^1$ is a
point on the unit circle. It readily follows that
\[
\va(r,\vs)=\left[\dfrac{8R^2 + 4R\varphi_1(r) s_2+ \varphi_1^2}{4R\sqrt{\varphi_1^2+ 4R^2+ 4R \varphi_1(r) s_2}}\right]^2, 
\qquad \vb(r,\vs)= \dfrac{2R^2\,\varphi_1(r)^2}{8R^2 + \varphi_1(r)^2+ 4R \varphi_1(r)\, s_2},
\]
\[
 \vl(r,\vs)= \dfrac{8\,R^2}{(8\,R^2+ \varphi_1^2+ 4R\varphi_1\,s_2)^2}, 
\qquad \textrm{and} \qquad(\va\,\vl)(r,\vs)= \dfrac{1}{2(4\,R^2 + \varphi_1^2 + 4R\,\varphi_1\,s_2)};
\]
hence in these new coordinates the energy surfaces $\Sigma_h$ can be written as 
\[
 \Sigma_h=\left\{(r,\vs,z_x,z_y) \in \R^+ \times \mathbb S^1 \times \R^2:\  \va(r,\vs) \, z_x^2+ z_y^2 = \dfrac{ r^2}{\vl(r,\vs)}
\Big(h - \Gamma/(8\pi) \, \log(\vb(r,\vs)\Big)\right\}.
\]
We also observe that 
\begin{itemize}
\item $\lim_{r\to 0^+} \va(r,\vs)=1\qquad \textrm{uniformly with respect to}\ \vs;$
\item $\lim_{r\to 0^+} \vb(r,\vs)=0\qquad \textrm{uniformly with respect to}\ \vs; $
\item $\lim_{r\to 0^+} \vl(r,\vs)=1/(8R^2)\qquad \textrm{uniformly with respect to}\ \vs. $
\end{itemize}

By taking into account the definition of the functions $\varphi_j$,
the right hand side of the equation defining the level set
$\Sigma_{h}$ reduces to
\[
\hat E(h,r, \vs):= \dfrac{\, r^2}{l(r,\vs)}\left[h- \Gamma/(8 \pi)\log(2 R^2 r^2 e^{-2/r^2}) + \Gamma/(8\pi)\log(\vc(r,\vs))\right]
\]
where $\vc(r,\vs):= 8R^2 + \varphi_1(r)^2+ 4R \varphi_1(r)\, s_2=8R^{2}+r^{2}e^{-2/r^{2}}+4Rre^{-1/r}s_{2}$.
We observe that 
\begin{equation}
\label{eq:energy_coll_man} \lim_{r \to 0^+} \hat E(h,r, \vs)=  \dfrac{2\,\Gamma\, R^2}{\pi},
\end{equation}
thus, as already implicit in Lemma \ref{thm:zonedinamiche}, in the
attractive case ($\Gamma >0$) the vortex point lies in the allowed region 
of every energy level $h$, while the opposite holds in the repelling case
($\Gamma<0$). However, a first important consequence of the change of
variable above introduced is that in the variables $(r,s)$ the kinetic
energy remains bounded when a collision occurs.

From now on, we shall only consider the attractive case. Thus we
assume
$$
\Gamma>0.
$$ 

The intersection between one (and hence every) energy hypersurface
$\Sigma_h$ with $r=0$ is called {\em total collision manifold.\/} In
virtue of the limit \eqref{eq:energy_coll_man}, we may
conclude that
\begin{itemize}
\item the total collision manifold does not depend on the fixed energy
  level $h$; otherwise stated it is a boundary of every energy level;
\item it is diffeomorphic to the two dimensional  torus $\mathbb T:=\mathbb S^1 \times \mathbb S^1$.
\end{itemize}

\begin{figure}
\subfigure[]{
\includegraphics[width=0.35\textwidth]{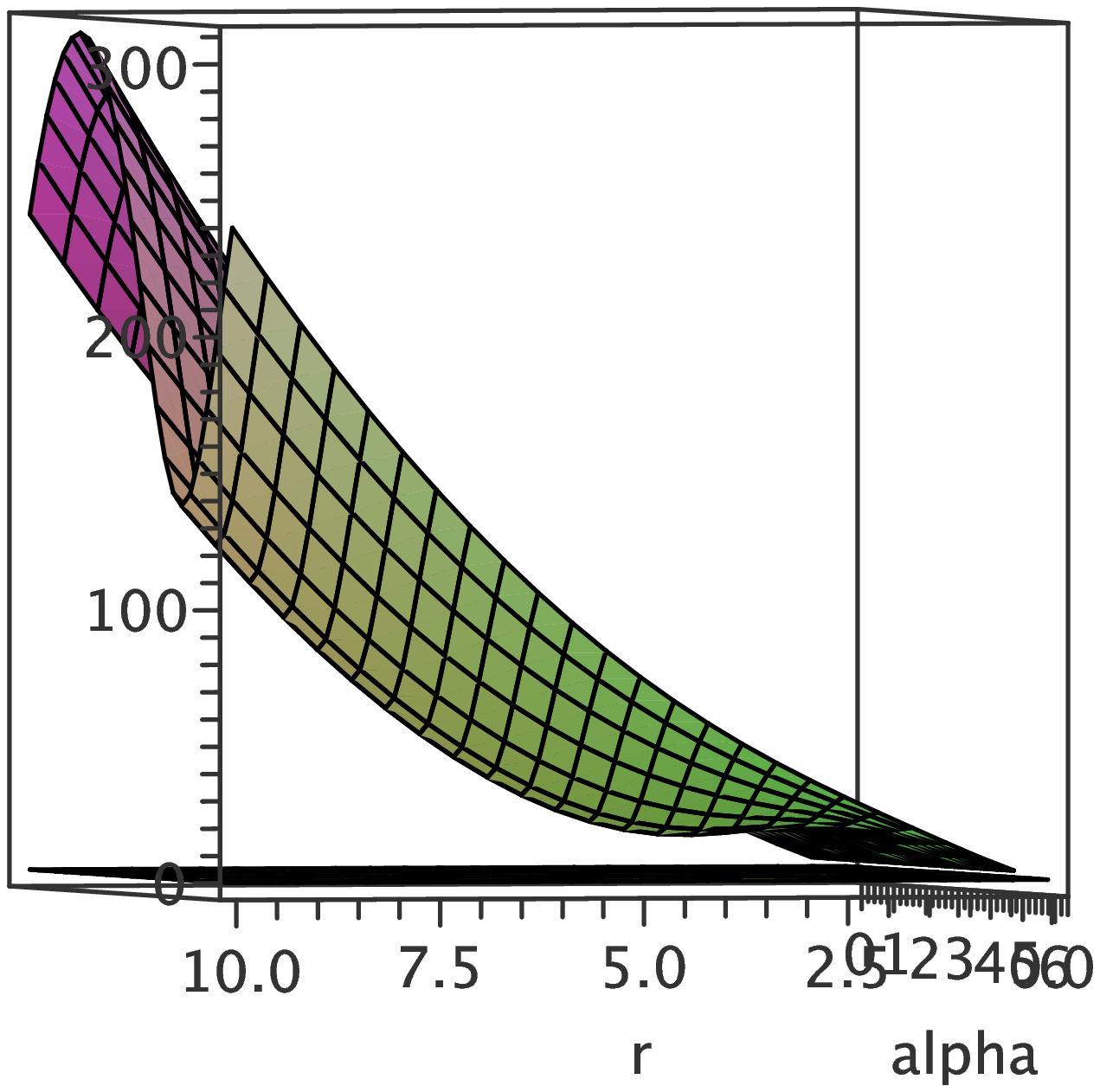}
}
\hspace{-28pt}
\subfigure[]{
\includegraphics[width=0.35\textwidth]{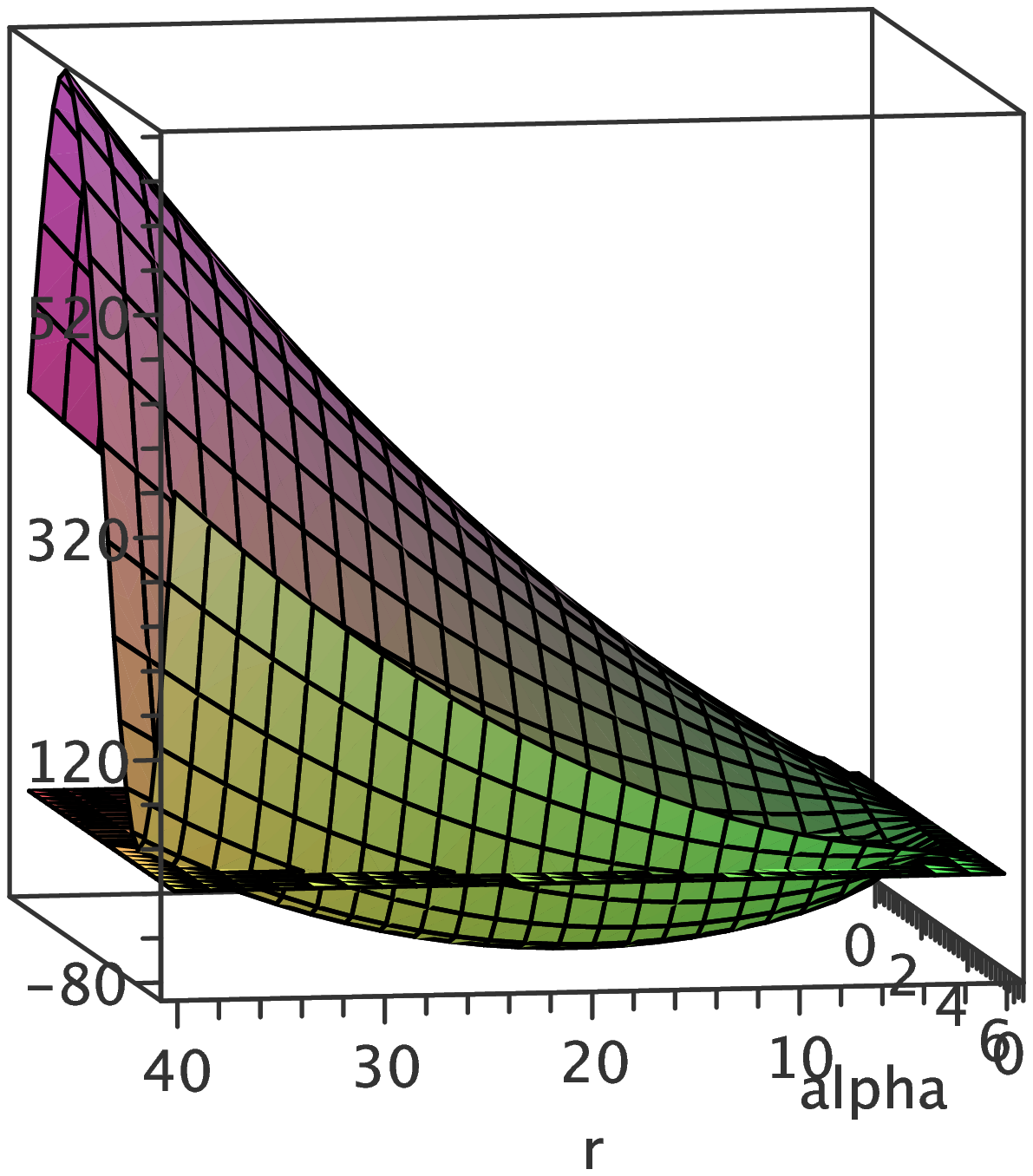}
}
\hspace{-28pt}
\subfigure[]{
\includegraphics[width=0.35\textwidth]{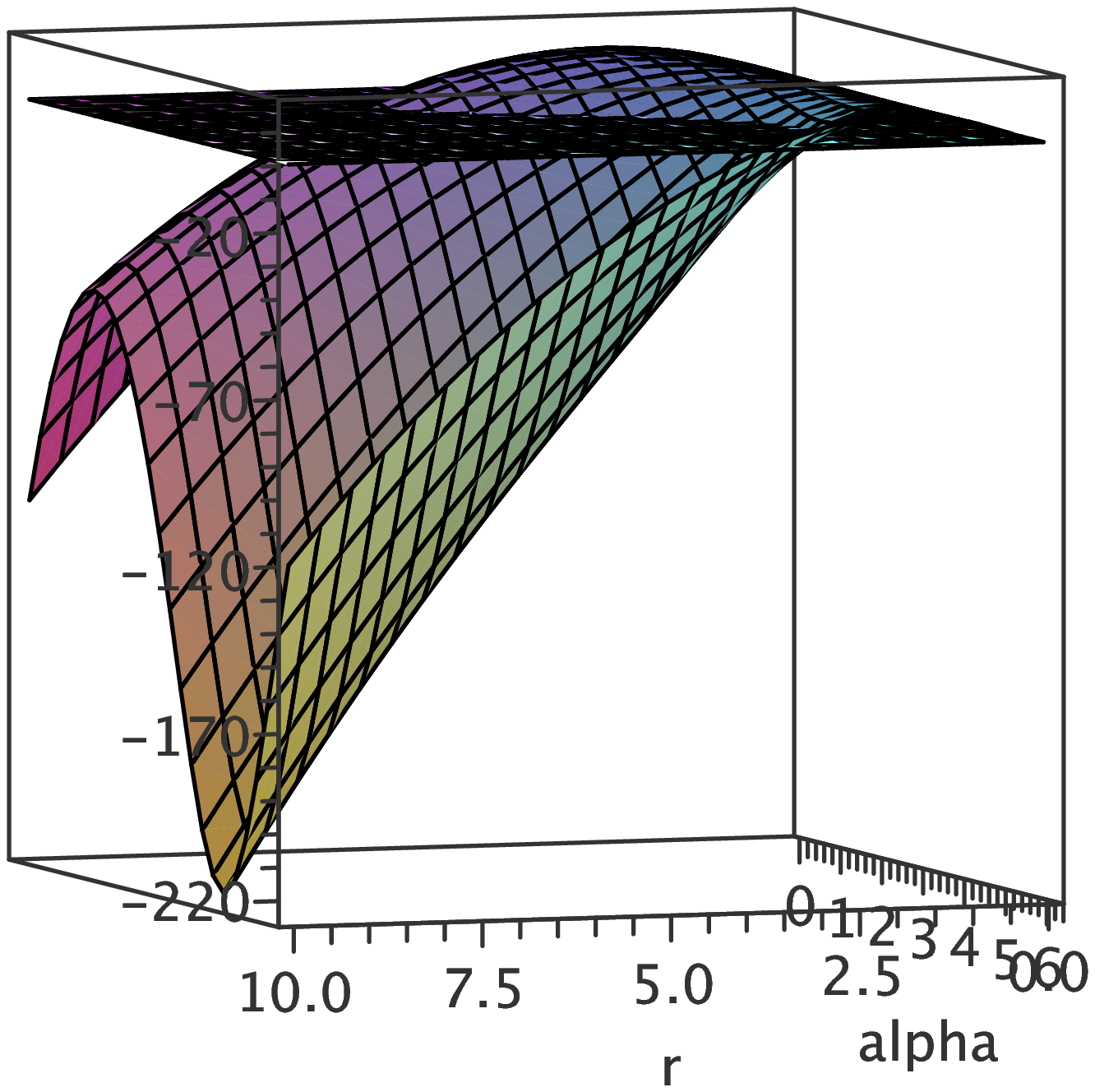}
}
\caption{Graph of the function $\hat E$ and of $\hat E=0$  in the three cases: (a) $h>h_{2}$, (b) $h_{1}\leq h \leq h_{2}$, (c) $h<h_{1}$.}
\label{fig:3casi}
\end{figure}

\begin{figure}[ht]
\centering
\subfigure[]
{
\includegraphics[width=0.300\textwidth]{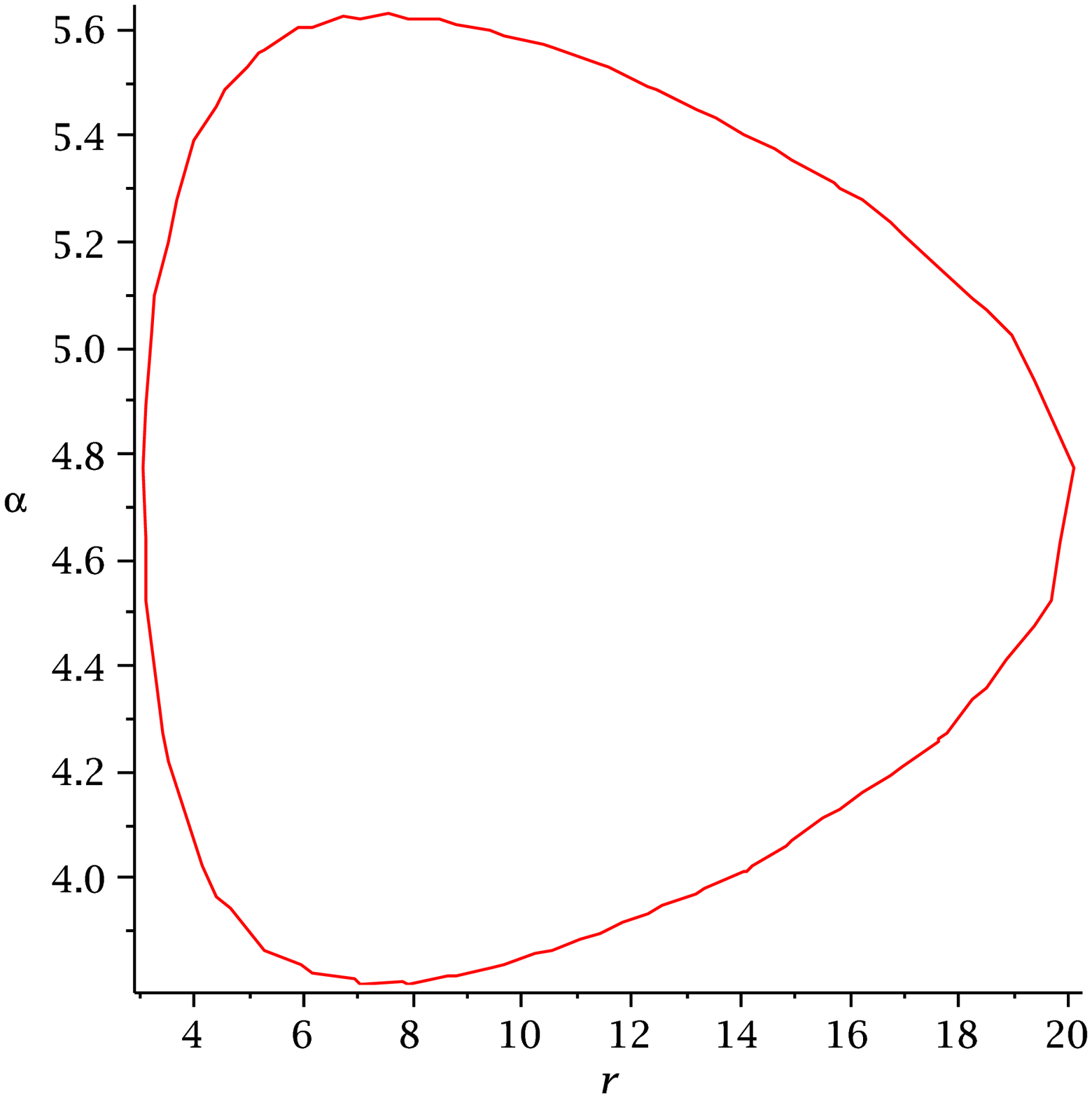}
}
\hspace{60pt}
\subfigure[]{
\includegraphics[width=0.300\textwidth]{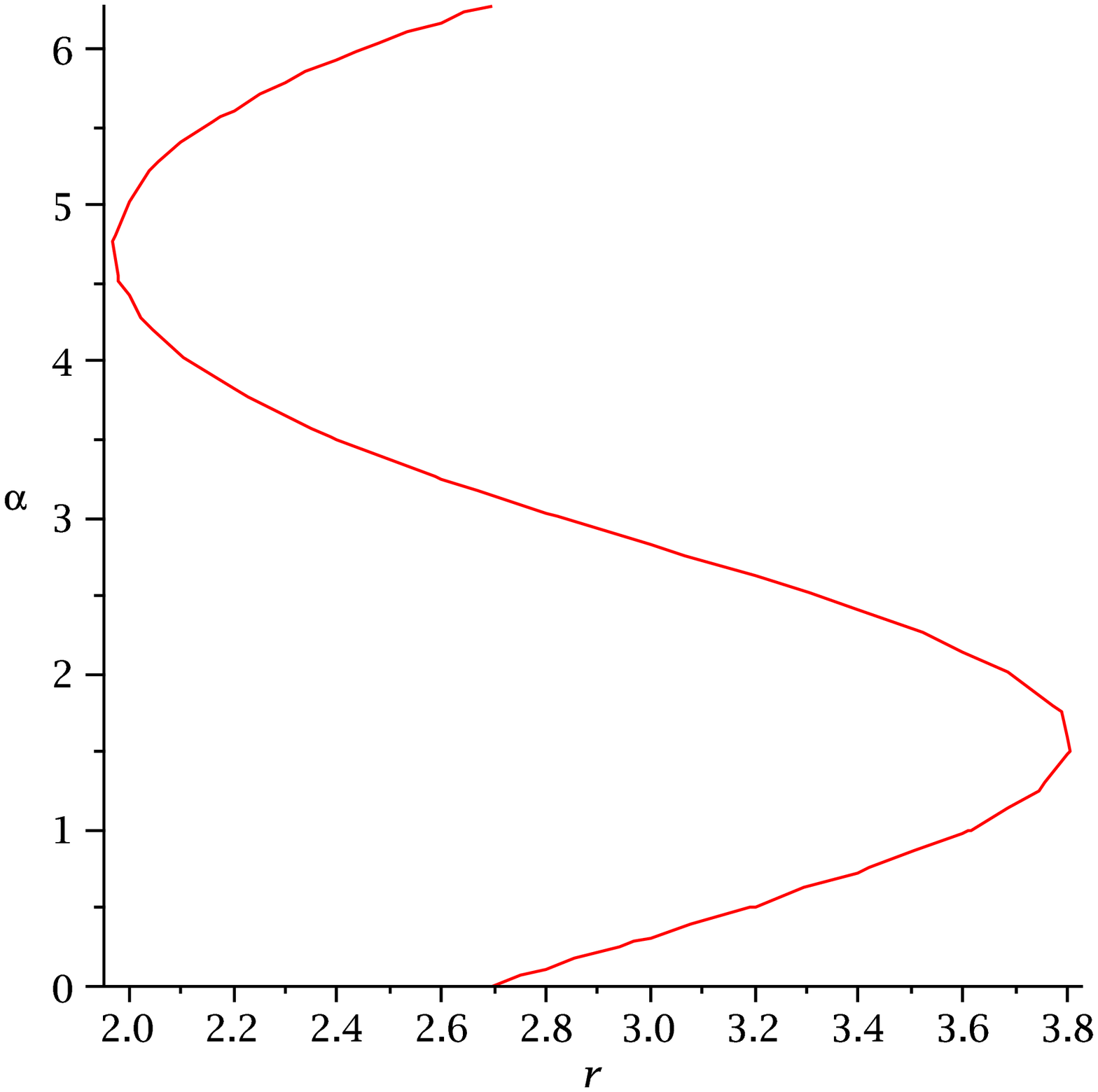}
}
\caption{Zero Velocity manifold (a) in the second case $h \in (h_1, h_2)$, (b) in the third case: $h<h_1$.}
\label{fig:zerovelnegative}
\end{figure}

From the dynamical viewpoint an important role is played by the zero
set of the function $\hat E$:
\[
\vZ_h:=\{(r,\alpha) \in X: \ \hat E(h,r,\alpha)=0\}.
\]
In the following we refer to this set as the {\em zero velocity
  manifold\/} in $\Sigma_h$. Rephrasing the results of Lemma
\ref{thm:zonedinamiche} in terms of the new coordinates $(r,\alpha)$
it readily follows that $\vZ_h$ is empty in the first case $(h>h_{2})$
and non-empty otherwise.  (Figure \ref{fig:3casi}). In the second
case, $(h_{1}\leq h \leq h_{2})$, the zero velocity manifold is
represented by a simple closed curve homeomorphic to a circle (or to a
point in the limit $h\to h_2$). The motion is forbidden in the region
bounded by the curve.  In the third case, the zero set can be seen as
the graph of a single-valued function $\alpha \mapsto r(\alpha)$ and the function
$\hat E$ is positive for $0<r<r(\alpha)$, which is the region on the left of the 
curve shown in figure (\ref{fig:zerovelnegative}b)

\subsection*{Regularization and McGehee coordinates}
We now use the new variables $r,$ $\alpha$ and $\vz$ in the equations
of motion \eqref{eq:hamiltonequ}. In order to preserve the continuity
of the flow with respect to the initial data, we need to ensure that
the transformed system has an everywhere differentiable vector
field. To this purpose we rescale the time variable in terms of the
distance from the singularity with the effect to exponentially
decrease the velocities near the singularity.  As a result the
collision solutions (which are singular in the old coordinates) move
along smooth orbits that asymptotically converge to the collision
manifold.

Let us define  $d\tau= \varphi_2(r)\varphi_1^{-1}(r)\, dt$ and use the notation 
\[
\langle \vz, \vs(\alpha)\rangle_a:= \va(x,y) z_x\, \cos\alpha + z_y\, \sin\alpha. 
\]
With the help of the identities 
\[
\dfrac{\varphi_1(r)}{\varphi_1'(r)}= \dfrac{r^3}{r^2+2}, \qquad  \dfrac{\varphi_1(r)}{\varphi_2^2(r)}= r^3\, e^{-1/r^2}, \qquad 
\dfrac{\varphi_1(r)\varphi_2'(r)}{\varphi_2(r)\varphi_1'(r)}=-\dfrac{r^2}{r^2+2}
\]
the Hamiltonian equations in \eqref{eq:hamiltonequ} become
\begin{equation}\label{eq:mcgehee1}
\left\{ \begin{array}{ll}
\dfrac{dr}{d\tau} &= \dfrac{2r^3}{(2+r^2)}\,\vl(r,\alpha)\,\langle \vz, \vs(\alpha)\rangle_a\\
\\
\dfrac{d\alpha}{d\tau}&= \vl(r,\alpha)\,(z_y \cos\alpha- \va(r, \alpha) z_x \sin \alpha)\\
\\
\dfrac{dz_x}{d\tau}&=- \, r e^{-1/r^2} [(\va\vl)_x(r, \alpha) z_x^2 + \vl_x(r,\alpha) \, z_y^2] - \dfrac{\Gamma}{8\pi}\, r^3\, e^{-1/r^2}
\dfrac{\vb_x(r,\alpha)}{\vb(r,\alpha)}+\\
 &+2\,\dfrac{r^2}{r^2+2} \vl(r,\alpha)\langle \vz, \vs(\alpha)\rangle_a\, z_x\\
\\
\dfrac{dz_y}{d\tau}&=- \, r e^{-1/r^2} [(\va\vl)_y(r, \alpha) z_x^2 + \vl_y(r,\alpha) \, z_y^2] - \dfrac{\Gamma}{8\pi}\, r^3\, e^{-1/r^2}
\dfrac{\vb_y(r,\alpha)}{\vb(r,\alpha)}+\\
 &+2\,\dfrac{r^2}{r^2+2} \vl(r,\alpha)\langle \vz, \vs(\alpha)\rangle_a\, z_y
\end{array}\right.
\end{equation}
where the subscripts in $(\va\vl)_x(r,\vs), (\va\vl)_y(r,\vs)$
(resp. $\vb_x(r,\vs), \vb_y(r,\vs)$) denote the partial derivative
with respect to the old cartesian variables. The derivate function is
then evaluated in the new coordinates at the point $(r,\alpha)$. 
The equations above are no longer singular at $r=0$: in fact, by computing
${\displaystyle \frac{\vb_{x}}{\vb}(r,\alpha)}$ one finds
${\displaystyle \frac{\vb_{x}}{\vb}(r,\alpha)\sim\varphi_{1}^{-1}}$
as $r\to 0$. Thus the time change produces the effect to regularize
the singularity. In addition the $\{r=0\}$ manifold results to be
invariant.

From a naive point of view, the study of the flow on the collision
manifold could appear meaningless, since the manifold is the image of
just a singular point where the orbits cease to exists. In reality,
the properties of the flow on such manifold yield informations on
the behavior of the orbits {\em close to} the singularity.

In order to simplify system \eqref{eq:mcgehee1}, we introduce a
further change of coordinates and time rescaling. Using the energy
relation
\begin{equation}\label{eq:energyrel}
\va(r,\alpha)\, z_x^2 + z_y^2 = \hat E(h,r,\alpha), 
\end{equation}
let us define $\psi$ and $\sigma$ such that
\begin{equation}
\left\{
\begin{array}{ll}
z_x= \sqrt{\hat E(h ,r , \alpha)/\va(r, \alpha)}\, \cos \psi, \\
z_y= \sqrt{\hat E(h, r, \alpha)} \, \sin \psi \\
\noalign{\vskip 3pt}d\tau=\sqrt{\hat E(h, r,\alpha)}\, d\sigma
\end{array}\right.
\end{equation}
Let us denote with $A_1(r,\alpha,z)$, $A_2(r,\alpha, z)$ the right
hand side of, respectively, the third and fourth equation in
\eqref{eq:mcgehee1}, so that $\dfrac{dz_z}{d\tau}=A_{1},
\dfrac{dz_y}{d\tau}=A_{2}$.
Therefore, on any fixed-energy shell, the system given in \eqref{eq:mcgehee1} 
reads 
\begin{equation}\label{eq:mcgehee2}
\left\{ \begin{array}{ll}
\dfrac{dr}{d\sigma} = \dfrac{2r^3\, \hat E(h, r, \alpha)}{(2+r^2)}\,\vl(r,\alpha)(\sqrt{\va(r,\alpha) }\, \cos\psi\, \cos \alpha\,+ \, \sin\psi\,\sin\alpha)\\
\\
\dfrac{d\alpha}{d\sigma}= \hat E(h,r,\alpha)\,\vl(r,\alpha)(\sin\psi\, \cos\alpha- \sqrt{\va(r, \alpha)}\, \cos\psi \sin \alpha)\\
\\\
\dfrac{d\psi}{d\sigma}=  B(r, \alpha, \psi)
\end{array}\right.
\end{equation}
where $B$ is given by:
\[
B(r,\alpha,\psi):= -\sqrt{\va} \, \sin\psi\, A_1+ \sqrt{\va}\sin \psi\, \cos\psi \dfrac{d}{d\tau}\,
\left(\sqrt{\dfrac{\hat E}{\va}}\right) + A_2\, \cos\psi - \dfrac{d}{d\tau}(\sqrt{\hat E}) \,\sin \psi \cos\psi.
\]


\subsection*{Flow and  invariant manifolds}

Recalling that $\vl$ is everywhere positive, the restpoints of
\eqref{eq:mcgehee2} correspond to solutions of the
following systems:
\begin{equation}\label{systems}
\left\{
\begin{array}{ll}
r=0\\
f_2(r,\alpha, \psi)=0\\
B(r,\alpha,\psi)=0
\end{array}\right. \qquad \textrm{or} \qquad
\left\{
\begin{array}{ll}
f_1(r,\alpha, \psi)=0\\
f_2(r,\alpha, \psi)=0\\
B(r,\alpha,\psi)=0
\end{array}\right. \qquad \textrm{or}\qquad  
\left\{
\begin{array}{ll}
\hat E(h,r,\alpha)=0\\
B(r,\alpha,\psi)=0
\end{array}\right.
\end{equation}
where 
\begin{equation*}
f_1(r,\alpha, \psi):= \sqrt{\va} \, \cos\psi\, \cos \alpha\,+ \, \sin\psi\,\sin\alpha, \qquad
f_2(r,\alpha, \psi) := \sin\psi\, \cos\alpha- \sqrt{\va}\, \cos\psi \sin \alpha\ .
\end{equation*}
We immediately discard the second system as it allows no solutions. 
We note that $\va \to 1$ for  $r\to 0$, thus the first system reduces to
\begin{equation*}
\left\{
\begin{array}{ll}
r=0\\
\sin(\psi-\alpha)=0\\
B(r,\alpha,\psi)=0
\end{array}\right.
\end{equation*}
whose solutions correspond to fixed points on the collision
manifold. The existence of solutions of the last system depends on the
energy level $h$: if $h\geq h_{2}$ the zero set of $\hat E$ is empty
and no solutions exist. For $h\leq h_2$, some solution
may exist.
 
Summarizing, any restpoint either lies on the collision manifold or on
the zero velocity manifold.  Let us first consider the collision
manifold: the asymptotic analysis of the function $B(r,\alpha, \psi)$
on the collision manifold (see Appendix A) gives
\[
\lim_{r\to 0^+} B(r, \alpha, \psi)=0.
\]
It follows that
\begin{lem}\label{thm:restpoints}
  The equilibria of the vector field given in \eqref{eq:mcgehee2} lying
  on the total collision manifold consists of two curves. In local
  coordinates $(r, \alpha, \psi)$ these curves are given by
\begin{enumerate}
\item[(i)]  $\mathscr P_1\equiv(0,\alpha,\alpha)$;
\item[(ii)]   $\mathscr P_2 \equiv (0,\alpha, \pi+\alpha)$.
\end{enumerate}
\end{lem}
\begin{prop}\label{pro:piep2selledegeneri}
  For each $\alpha$, the equilibrium points
\[
(0,\alpha,\alpha) \in \mathscr P_1
\]
and 
\[
(0, \alpha, \pi+\alpha) \in \mathscr P_2 
\]
are degenerate saddles. 
\begin{enumerate}
\item $
 \dim W^u(\mathscr P_1)=1, \qquad  \dim W^s(\mathscr P_1)=1, \qquad  \dim W^0(\mathscr P_1)=1.$
\item $
\dim W^u(\mathscr P_2)=1, \qquad  \dim W^s(\mathscr P_2)=1, \qquad  \dim W^0(\mathscr P_2)=1.$
\end{enumerate}
\end{prop}
\proof The flow on the collision manifold is given by
\begin{equation}\label{eq:mcgehee22trisdiel=2collmanfld}
\left\{
\begin{array}{ll}
\dfrac{d\alpha}{d\sigma}= \dfrac{\Gamma}{4\pi}\, \sin (\psi-\alpha)\\
\\
\dfrac{d\psi}{d \sigma} = 0.\\
\end{array}\right.
\end{equation}
whose orbits are parallel to $\alpha$-axis and flow from
$\mathscr P_{2}$ to $\mathscr P_{1}$. The stability of the restpoints
is determined by the eigenvalues of the Jacobian matrix of
\eqref{eq:mcgehee2}. It follows (see appendix \ref{subsec:asy}) that
for any point $P_{1}\in\mathscr P_{1}$ and $P_{2}\in\mathscr P_{2}$
the eigenvalues are
 \begin{equation*}
P_{1}\in\mathscr P_{1}\Rightarrow\left\{
\begin{array}{ll}
\lambda_r=0\\
\lambda_\alpha=-\dfrac{\Gamma}{4\pi}\\
\lambda_\psi=0. 
\end{array}\right.,\qquad
P_{2}\in\mathscr P_{2}\Rightarrow \left\{
\begin{array}{ll}
\lambda_r=0\\
\lambda_\alpha=\dfrac{\Gamma}{4\pi}\\
\lambda_\psi=0. 
\end{array}\right.
\end{equation*}

\begin{figure}[htbp]
\centering
\includegraphics[width=0.45\textwidth]{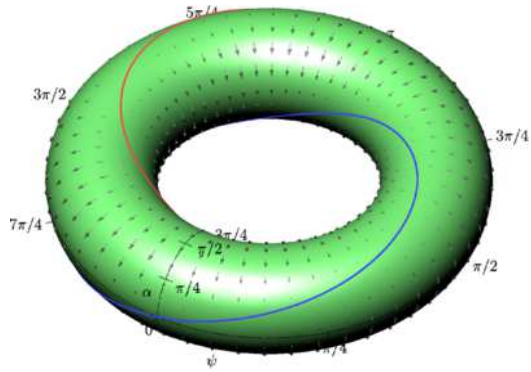}
\caption{The collision manifold, the curves of restpoints  $\mathscr P_{1}$ (blue) and $\mathscr P_{2}$ (red), and the vectorfield of equations \eqref{eq:mcgehee22trisdiel=2collmanfld}.}
\label{fig:streamlines}
\end{figure}

and are coherent with the dynamics restricted on the collision
manifold as given by \eqref{eq:mcgehee22trisdiel=2collmanfld}, where
$\mathscr P_{1}$ is an attractor and $\mathscr P_{2}$ is a repeller
(figure \ref{fig:streamlines}). However, the presence of null
eigenvalues implies that the linear approximation of the flow, taken
alone, does not provide enough information to determine the
qualitative dynamics close to the equilibrium points.  As flow in the
$\psi$ direction is null (and in fact $\psi$ can be regarded as a
parameter for an equilibrium point), in order to determine the
asymptotic behavior close to $P_1$ and $P_2$, it is enough investigate
the dynamics restricted to the $(r,\alpha)$-plane. The proof that the
two equilibrium curves are indeed degenerate saddles follows by direct
integration of the system once the equations have been expanded around
the equilibrium point in Taylor series. We omit the details and we
refer to the equivalent proof of Lemma 7.4 in \cite{stoica}.  \finedim

\begin{table}\centering
\begin{tabular}{|l|c|c|c|}
\hline & 
\begin{sideways}$\dim  W^{s}\ $\end{sideways} &
\begin{sideways}$\dim  W^{u}\ $\end{sideways} &
\begin{sideways}$\dim W^0\  $\end{sideways}
\\ \hline 
&&&\\
\multirow{1}{*}{At $\mathscr P_1$}& $1$  &
$1$ & $1 $  \\ 
\hline  &&&\\
\multirow{1}{*}{At  $\mathscr P_2$}  & $1$&
$1$ & $1$ \\ 
\hline
\end{tabular}
\caption{Dimensions of the invariant  manifolds along the equilibrium curves $\mathscr P_1$ and $\mathscr P_2$}
\label{tb:tabella1}
\end{table}

\begin{defn}
  We shall say that the flow on the collision manifold is {\em totally
    degenerate\/} if the unstable manifold of an equilibrium point
  $\vP_1 \in \mathscr P_1$ , coincides with the stable manifold of
  some equilibrium point $\vP_2 \in \mathscr P_2$.
\end{defn}
\begin{lem}\label{thm:flussodegcollmnfld}
  The flow on the total collision manifold is totally degenerate. More
  precisely
\begin{enumerate}
\item[(i)] $W^u(\vP_1) \equiv W^s(\vP_2)$; 
\item[(ii)] $W^u(\vP_2) \equiv W^s(\vP_1)$;
\end{enumerate}
where $\vP_1 \in \mathscr P_1$ and $\vP_2 \in \mathscr P_2$ are chosen
in such a way the last coordinate of the two points agrees.
\end{lem}
\proof The proof of this result follows by a straightforward
integration of the equations of motion on the total collision manifold
$r=0$.
\finedim\\
A direct consequence of the previous result is the following:
\begin{cor}({\bf Existence of heteroclinic
    connections\/})\label{thm:collhetoncollmanfold} There exists an
  heteroclinic connection between each equilibrium point $\vP_1 \in
  \mathscr P_1$ and the point $\vP_2 \in \mathscr P_2$ where $\vP_1,
  \vP_2$ were chosen in such a way that they have the same projection
  on the first and third coordinate.
\end{cor}
\proof The proof of this result follows immediately by the previous
result. By the fact that $r=0$ and $\psi$ is constant, it follows that
the non equilibrium solutions are in the $(\alpha, \psi)$-plane lines
parallel to the $\alpha$-axis.  Moreover each point of equilibrium on
$\mathscr P_1$ is attracting while each equilibrium point on $\mathscr
P_2$ is repelling.
\finedim \\

Moving out of the collision manifold, the two lines $\mathscr P_{1}$
and $\mathscr P_{2}$ exhibit the opposite stability character: indeed
$\dfrac{dr}{d\sigma}>0$ for $\psi=\alpha$ and $\dfrac{dr}{d\sigma}<0$
when $\psi=\alpha+\pi$, meaning that the system goes into the
collision along $\mathscr P_{2}$ and escape from the collision along
$\mathscr P_{1}$.

Next we examine the restpoints and the flow on the zero velocity
manifold. This is more easily accomplished by looking at the system
given in \eqref{eq:mcgehee1}. Restpoints, in fact, are not changed by
a time scaling. Since the zero velocity manifold coincides with the
zero set of the function $\hat E$ and by taking into account 

Setting $\hat E=0$ in the energy relation \eqref{eq:energyrel}, which
implies $z_x=z_y=0$, it follows that on the zero velocity manifold the
dynamical system \eqref{eq:mcgehee1} reduces to:
\begin{equation}\label{eq:mcgeheezerovelmnfld}
\left\{ \begin{array}{ll}
\dfrac{dr}{d\tau} = 0\\
\\
\dfrac{d\alpha}{d\tau}= 0\\
\\
\dfrac{dz_x}{d\tau}=- \dfrac{\Gamma}{8\pi}\, r^3\, e^{-1/r^2}
\dfrac{\vb_x(r,\alpha)}{\vb(r,\alpha)}\\
\\
\dfrac{d z_y}{d\tau}= - \dfrac{\Gamma}{8\pi}\, r^3 e^{-1/r^2}
\dfrac{\vb_y(r,\alpha)}{\vb(r,\alpha)}.
\end{array}\right.
\end{equation}
The restpoints on the zero velocity manifold (if any) correspond 
to the solutions of the equations:
\[
\dfrac{ \vb_x(r,\alpha)}{\vb(r,\alpha)}=\dfrac{\vb_y(r,\alpha)}{\vb(r,\alpha)}=0.
\]
An elementary calculation shows the following result.
\begin{lem}\label{thm:restvelmnfld}
  For $h=h_2$ (defined above as $h_2:= \dfrac{\Gamma}{4\pi} \log
  (2R)$) there exists only one restpoint on the zero velocity
  manifolds at $P:=(r_*, 3\pi/2)$ for $\varphi_1(r_*)= 4R$.  For
  $h\neq h_2$ there are no restpoints.
\end{lem}

\section{Global flow and dynamics on the sphere}\label{sec:global}

It is now possible to bring back on the sphere the results found on
the stereographic plane in the previous sections.

In terms of the coordinates $(\phi,\theta)$ on the sphere the Hamiltonian reads as:
\[
 H(\theta, \phi, p_\theta, p_\phi)= \dfrac{1}{2 R^2} \left(\dfrac{1}{\sin^2\theta} \, p_\phi^2+ p_\theta^2\right) + \dfrac{\Gamma}{8\pi} \log(2R^2(1-\sin\theta\sin\phi)).
\]
Note that 
the vortex is located at $(\phi,\theta)=(\pi/2, \pi/2)$,
therefore for $(\phi,\theta)\to (\pi/2, \pi/2)$ the dynamical behavior
becomes unknown since the vectorfield ceases to exist.

Let us call {\em vortex half-sphere\/} the half-sphere centered around
the vortex point, and {\em antivortex half-sphere\/} the complementary
half-sphere; let us call {\em vortex-parallel\/} any circle on the
sphere equidistant from the vortex and {\em vortex-meridian\/} any
great circle passing through the vortex. Finally let us call {\em
  antipodal point} the point on the sphere opposite to the vortex
point.  On the sphere the results of Lemma \ref{thm:zonedinamiche} can
be rephrased as follows.
 \begin{thm}\label{thm:1}
   If $h<h_{2}$ the motion is allowed in the region of the sphere
   containing the vortex point and bounded by a vortex-parallel that
   lies on the vortex half-sphere for $h<h_{1}$ and in the
   antivortex half-sphere otherwise.  If $h\geq h_{2}$ then the motion
   is allowed everywhere on the
   sphere.
\end{thm}
Moreover lemma \ref{thm:restvelmnfld} is rephrased as:
\begin{thm}\label{thm:2}
  For $h=h_{2}$ the zero velocity manifold consists only of one
  point which is the antipodal point.
\end{thm}

In order to understand the global dynamics it is useful to show the
existence of a second conserved quantity, analogous to the angular
momentum for planar dynamics.  To this aim it is convenient to move
the vortex point at the north pole $N=(0,0,2R)$ (or, equivalently, to
redefine the parameterization of the sphere). Obviously this does not
change the dynamics. Note that the curves $\{\theta=const\}$ and
$\{\phi=const\}$ now respectively correspond to the vortex parallels
and to the vortex meridians.

In this setting, the dynamical system reads as
\begin{equation}\label{eq:HJonspherenord}
\left\{
\begin{array}{ll}
\dot \phi = \dfrac{1}{R^2\, \sin^2\theta} \, p_\phi\\
\\
\dot \theta= \dfrac{p_\theta}{R^2}\\
\\
\dot p_\phi=0\\
\\
\dot p_\theta=\dfrac{\cos\theta}{R^2\, \sin^3\theta} \,p_\phi^2+\dfrac{\Gamma}{8\pi}\dfrac{\sin\theta}{2R^{2}(1+\cos\theta)}\\
\\
\end{array}\right. 
\end{equation}
corresponding to the Hamiltonian
$$
H(\theta, \phi, p_\theta, p_\phi)= \dfrac{1}{2 R^2} \left(\dfrac{1}{\sin^2\theta} \, p_\phi^2+ p_\theta^2\right) + \dfrac{\Gamma}{8\pi} \log(2R^2(1+\cos\theta)).
$$
Going to the Lagrangian formulation, we can write the
Euler-Lagrange equations
$$
\frac{d}{dt}(R^{2}\dot\phi\sin^{2}\theta)=0
$$  
$$
R^{2}\ddot\theta-R^{2}\sin\theta\cos\theta\dot\phi^{2}+\frac{\Gamma}{8\pi}\frac{\sin\theta}{2R^{2}(1+\cos\theta)}=0.
$$
The first shows the existence of a conserved quantity, namely the
spherical angular-momentum $l=R^{2}\sin^{2}\theta\dot\phi$. It follows
that
\begin{lem}
  A necessary condition for a solution to either collide with the
  vortex or to reach the antipodal point is $l=0$.
\end{lem}

\proof Writing the energy relation in terms of
$(\phi,\theta,\dot\phi,\dot\theta)$ and substituting
$\dot\phi=\frac{l}{R^{2}\sin^{2}\theta}$, it follows that a solution
exists only for those $\theta$ satisfying
$$
2R^{2}h\sin^{2}\theta-\frac{\Gamma}{4\pi}R^{2}\sin^{2}\theta\log(2R^{2}(1+\cos\theta))-l^{2}\geq 0
$$
Recalling that the vortex is placed at $\theta=\pi$, and that the
antipodal point is at $\theta=0$, it follows that if either a
collision occurs, or the point goes to the antipodal point, then
$l^{2}\leq 0$.
\finedim\\

\begin{figure}
\centering
{
\includegraphics[width=0.840\textwidth]{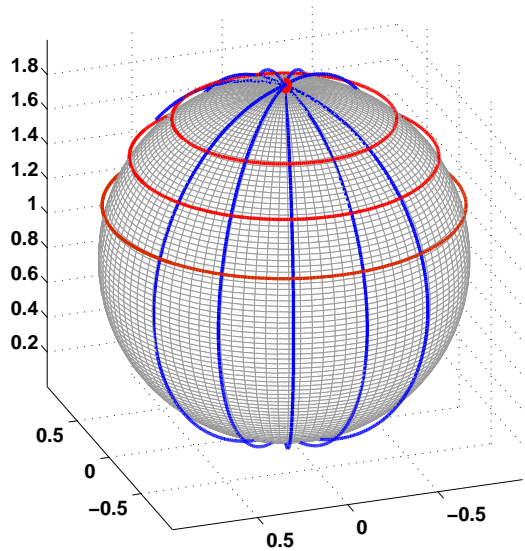}
}
\caption{Particular orbits lying on the vortex-paralles and vortex-meridians }
\label{fig:orbitesfera}
\end{figure}

Looking at the system \eqref{eq:HJonspherenord}, one can easily prove
the existence of particular solutions as depicted in
Fig.\ref{fig:orbitesfera}

\begin{lem}
\begin{itemize}
\item[(i)] Any vortex-parallel on the vortex half-sphere is the
  support of a periodic orbit.
\item[(ii)] The vortex-meridians are flow-invariant.
\end{itemize}
\end{lem}
\proof
For any $\simbolovettore\theta\in(\tfrac{\pi}{2},\pi)$ the curve 
\[
\gamma_{\simbolovettore{\theta}}(t):=(\phi,\theta,p_{\phi},p_{\theta})(t)=\left(\phi_{0}+t\frac{\simbolovettore{p_{\phi}}}{R^{2}(\sin\simbolovettore\theta)^{2}},
\simbolovettore\theta,\simbolovettore{p_{\phi}},0\right)
\]
with
$\simbolovettore{p_{\phi}}^{2}=-\frac{\Gamma}{8\pi}\frac{(\sin\simbolovettore{\theta})^{4}}
{2\cos\simbolovettore{\theta}(1+\cos\simbolovettore{\theta})}$ is a
solution of the system. Note that the previous relation can not be
satisfied if $\simbolovettore{\theta}\in(0,\pi/2]$, which implies that only
the vortex-parallel placed in the vortex half-sphere are support of
periodic orbits.  Moreover the period of
$\gamma_{\simbolovettore{\theta}}(t)$ tends to zero as
$\simbolovettore{\theta}$ goes to $\pi/2$ or $\pi$. This proves statement
$(i)$.  Statement $(ii)$ immediately follows by noting that
any initial data
$(\phi,\theta,p_{\phi},p_{\theta})(0)=(\phi^{0},\theta^{0},0,p_{\theta}^{0})$
leads to an orbit traveling on the $\{\phi=\phi_{0}\}$ vortex
meridian.
\finedim\\

A consequence of the previous lemmas is that any orbit with energy
$h>h_{2}$ that passes through the antipodal point will end into the
collision. 

The existence of heteroclinic connection anywhere on the total
collision manifold of the regularized flow (Corollary
\ref{thm:collhetoncollmanfold}) provides the way to extend beyond the
collision the orbits of the singular flow. In fact, if
$\gamma(t)=(\phi,\theta)(t):[0,T_{s})\rightarrow \vS$ is a collision
solution ending in the singularity at time $T_{s}$, we define the
collision-transmission solution as the path
$\bar\gamma:[0,2T_{s}]\rightarrow \vS$ as
$$
\bar\gamma(t):=\left\{\begin{array}{ll}
\gamma(t),&t\in[0,T_{s})\\
(\phi_{V},\theta_{V}) &t=T_{s}\\
(2\phi_{V}-\phi(2T_{s}-t),2\theta_{V}-\theta(2T_{s}-t)) &t\in(T_{s},2T_{s}]\\
\end{array}\right.
$$
where $(\phi_{V},\theta_{V})$ are the coordinates of the vortex
point. The extended flow obtained by replacing the singular
trajectories with the collision-transmission solution results to be
continuous with respect to the initial data. The same result for a
single logarithmic center on the plane has already been proved, among
others, in \cite{cate} with a completely different technique.

Finally, from the above discussion it follows that the
collision-transmission solution behaves in three different ways,
depending on the energy level $h$: if $h<h_{2}$, after the ejection
from the singularity, the particle reaches the zero velocity manifold,
then it reverses the motion and falls back into the vortex point; if
$h=h_{2}$ after the ejection the orbit asymptotically reaches the
antipodal restpoint (this is an heteroclinic orbit between a point of
$\mathscr{P}_2$ and the single restpoint on the zero velocity
manifold), and if $h>h_{2}$ after the ejection, the orbits travels
along a vortex meridian, passes through the antipodal restpoint,
continues the motion on the opposite meridian and falls down again
into the singularity.

\appendix
\section{Useful asymptotics } \label{sec:asympt}

In this appendix we list some asymptotic limits of the functions
appearing in the equation of motion in McGeheee coordinates. They are
useful to compute the spectrum of the eigenvalues associated to the
fixed points.

\subsection*{On the total collision manifold $r=0$}\label{subsec:asy}

All the limits below are computed with respect to $r$ and are
uniform with respect to the other variables.  For the function $\hat
E$, we have:
\[
\lim_{r \to 0^+} \hat E(h, r, \alpha)= \dfrac{4\Gamma\, R^4}{\pi}\quad
\lim_{r \to 0^+}  \hat \partial_r \,E(h, r, \alpha)=0\quad  \lim_{r \to 0^+} \hat \partial_\alpha \,E(h, r, \alpha)=0
\]
For the functions $\va, \vb$, we have 
\[
\lim_{r \to 0^+} \va( r, \alpha)= 1,\quad
\lim_{r \to 0^+} \partial_r\va( r, \alpha)=0,\quad
\lim_{r \to 0^+} \partial_\alpha \va( r, \alpha)=0
\]
\[
\lim_{r \to 0^+} \vb(r,\alpha)=0,\quad
\lim_{r \to 0^+}\partial_r \vb(r,\alpha)=0
\]
For the functions $\va_x, \vb_x, \va_y, \vb_y$, we have 
\[
\lim_{r \to 0^+} \va_x( r, \alpha)= 0\quad
\lim_{r \to 0^+} \va_y( r, \alpha)= 0\quad
\lim_{r \to 0^+} \partial_r\va_x( r, \alpha)=0
\]
\[
\lim_{r \to 0^+} \vb_x(r,\alpha)=0\quad
\lim_{r \to 0^+} \vb_y( r, \alpha)= 0\quad
\lim_{r \to 0^+}\partial_r \vb_x(r,\alpha)=0
\]
\[
\lim_{r\to 0^+}  \dfrac{r^3\, e^{-1/r^2} \vb_x(r,\alpha)}{\vb(r,\alpha)}=0
\]

\[
\lim_{r \to 0^+} \partial_\alpha\va_x( r, \alpha)=0\\
\lim_{r \to 0^+} \partial_\alpha\va_y( r, \alpha)=0
\]
\[
\lim_{r \to 0^+} \partial_\alpha\vb_x( r, \alpha)=0\\
\lim_{r \to 0^+} \partial_\alpha\vb_y( r, \alpha)=0
\]
\[
\lim_{r\to 0^+}  \dfrac{r^3\, e^{-1/r^2} \partial_\alpha \vb_x(r,\alpha)}{\vb(r,\alpha)}=0,\qquad  
\lim_{r\to 0^+}  \dfrac{r^3\, e^{-1/r^2} \partial_\alpha \vb_y(r,\alpha)}{\vb(r,\alpha)}=0
\]
\[
\lim_{r\to 0^+}  \dfrac{r^3\, e^{-1/r^2} \vb_x(r,\alpha)\partial_\alpha \vb_x(r,\alpha)}{\vb^2(r,\alpha)}=0,\qquad  
\lim_{r\to 0^+}  \dfrac{r^3\, e^{-1/r^2} \vb_y(r,\alpha)\partial_\alpha \vb_y(r,\alpha)}{\vb^2(r,\alpha)}=0
\]

For the functions $z_x, z_y$, we have:
\[
\lim_{r \to 0^+} z_x(r,\alpha, \psi)=\sqrt{\dfrac{\Gamma}{\pi}}\,2 R^2 \cos\psi, 
\quad \lim_{r\to 0^+} z_y(r,\alpha, \psi)= \sqrt{\dfrac{\Gamma}{\pi}}\,2 R^2 \sin\psi,
\]
\[
\lim_{r\to 0^+} \partial_r z_x(r,\alpha, \psi)=0,\quad \lim_{r \to 0^+} \partial_\alpha z_x(r,\alpha, \psi)=0, 
\]
\[
\lim_{r\to 0^+} \partial_\psi z_x(r,\alpha, \psi)=-\sqrt{\dfrac{\Gamma}{\pi}}\,2 R^2 \sin\psi,\quad \lim_{r \to 0^+} \partial_\psi z_y(r,\alpha, \psi)= 
\sqrt{\dfrac{\Gamma}{\pi}}\, 2R^2 \cos\psi, 
\]
\[
\lim_{r\to 0^+} \partial_r z_y(r, \alpha, \psi)=0, \quad \lim_{r \to 0^+} \partial_\alpha z_y(r,\alpha, \psi)=0. 
\]
As consequence of the above asymptotic behavior it follows that 
\[
\lim_{r \to 0^+} A_1(r, \alpha, \psi)=0, \qquad  \lim_{r \to 0^+} A_2(r, \alpha, \psi)=0, \qquad \lim_{r \to 0^+} B(r, \alpha, \psi)=0
\]
\[
\lim_{r \to 0^+} \partial_\alpha A_1(r, \alpha, \psi)=0, \qquad  \lim_{r \to 0^+} \partial_\alpha A_2(r, \alpha, \psi)=0.
\]
Denoting by $J:=(J_{ij})_{i,j}$ the variational matrix on the total
collision manifold it follows that
\[
J_{11}=0, \quad J_{12}=0, \qquad J_{13}=0, \qquad J_{32}=0, \qquad J_{33}=0.
\]
\[
J_{21}=0, \qquad J_{22}=-\dfrac{\Gamma}{2\pi}\, \cos(\psi-\alpha),\qquad J_{23}=-J_{22},
\]
\[
\lim_{r\to 0^+} \partial^2_{\alpha \tau} \va(r,\alpha)=0, \qquad \lim_{r\to 0^+} \partial^2_{\alpha \tau} \hat E(h,r,\alpha)=0.  
\]
The limit involved in the computation of the term $J_{31}$, in general
may not exists. However at the restpoints $\psi=\alpha$ or
$\psi=\alpha+\pi$ this limit actually exists and this implies that
$J_{31}=0$.

\end{document}